\documentclass[12pt]{amsart}
\usepackage{amscd,amssymb}
\usepackage{amsthm,amsmath,amssymb}
\usepackage[matrix,arrow]{xy}

\theoremstyle{definition}
\newtheorem{theorem}[equation]{Theorem}

\newtheorem{proposition}[equation]{Proposition}
\newtheorem{lemma}[equation]{Lemma}
\newtheorem{corollary}[equation]{Corollary}

\newtheorem{example}[equation]{Example}
\newtheorem{definition}[equation]{Definition}

\theoremstyle{remark}
\newtheorem{remark}[equation]{Remark}

\makeatletter\@addtoreset{equation}{section} \makeatother

\begin{document}

\title{Halphen pencils on quartic threefolds}\footnote{The work was partially supported by
RFFI grant No. 08-01-00395-a~and grant N.Sh.-1987.2008.1.}

\author{Ivan~Cheltsov and Ilya~Karzhemanov}

\thanks{We assume that all varieties
are projective, normal and  defined over $\mathbb{C}$.}

\pagestyle{headings}

\begin{abstract}
For any smooth quartic threefold in $\mathbb{P}^{4}$ we classify
pencils on it whose general element~is an~irreducible surface
birational to a~surface of Kodaira~dimension zero.
\end{abstract}

\maketitle

\section{Introduction}
\label{section:into}

Let $X$ be a~smooth quartic threefold in $\mathbb{P}^{4}$. The
following result is proved in \cite{IsMa71}.

\begin{theorem}
\label{theorem:IM} The threefold $X$ does not contain pencils
whose general element is an~irreducible surface that is birational
to a~smooth surface of Kodaira~dimension $-\infty$.
\end{theorem}

On the~other hand, one can easily see that the~threefold $X$
contains infinitely many pencils whose general elements are
irreducible surfaces of Kodaira dimension zero.

\begin{definition}
\label{definition:Halphen-pencil} A Halphen pencil is
a~one-dimensional linear system whose general element is
an~irreducible subvariety birational to a~smooth variety of
Kodaira~dimension~zero.
\end{definition}

The~following result is proved in \cite{ChPa05h}.

\begin{theorem}
\label{theorem:Halphen-pencils-quartic} Suppose  that $X$ is
general. Then every Halphen pencil on $X$ is cut out~by
$$
\lambda l_{1}\big(x,y,z,t,w\big)+\mu l_{2}\big(x,y,z,t,w\big)=0\subset\mathrm{Proj}\Big(\mathbb{C}[x,y,z,t,w]\Big)\cong\mathbb{P}^{4},%
$$
where $l_{1}$ and $l_{2}$ are linearly independent linear forms,
and $(\lambda:\mu)\in\mathbb{P}^{1}$.
\end{theorem}

The assertion of Theorem~\ref{theorem:Halphen-pencils-quartic} is
erroneously proved in \cite{Ch00a} without the~assumption that
the~threefold $X$ is general. On the~other hand, the~following
example is constructed~in~\cite{Isk01}.

\begin{example}
\label{example:Iskovskikh} Suppose that $X$ is given by the~
equation
$$
w^{3}x+w^{2}q_{2}\big(x,y,z,t\big)+wxp_{2}\big(x,y,z,t\big)+q_{4}\big(x,y,z,t\big)=0\subset\mathrm{Proj}\Big(\mathbb{C}[x,y,z,t,w]\Big)\cong\mathbb{P}^{4},%
$$
where $q_{i}$ and $p_{i}$ are forms of degree $i$. Let
$\mathcal{P}$ be the~pencil on $X$ that is cut out by
$$
\lambda x^{2}+\mu\Big(wx+q_{2}\big(x,y,z,t\big)\Big)=0,
$$
where $(\lambda:\mu)\in\mathbb{P}^{1}$. Then $\mathcal{P}$ is a~
Halphen pencil if $q_{2}(0,y,z,t)\ne 0$ by
\cite[Theorem~2.3]{ChPa05h}.
\end{example}

The purpose of this paper is to prove the~following result.

\begin{theorem}
\label{theorem:main} Let $\mathcal{M}$ be a Halphen pencil on $X$.
Then
\begin{itemize}
\item either $\mathcal{M}$ is cut out on $X$ by the~pencil
$$
\lambda l_{1}\big(x,y,z,t,w\big)+\mu l_{2}\big(x,y,z,t,w\big)=0\subset\mathrm{Proj}\Big(\mathbb{C}[x,y,z,t,w]\Big)\cong\mathbb{P}^{4},%
$$
where $l_{1}$ and $l_{2}$ are linearly independent linear forms, and $(\lambda:\mu)\in\mathbb{P}^{1}$, %

\item or the~threefold $X$ can~be given by the~equation
$$
w^{3}x+w^{2}q_{2}\big(x,y,z,t\big)+wxp_{2}\big(x,y,z,t\big)+q_{4}\big(x,y,z,t\big)=0\subset\mathrm{Proj}\Big(\mathbb{C}[x,y,z,t,w]\Big)\cong\mathbb{P}^{4}
$$
such that $q_{2}(0,y,z,t)\ne 0$, and $\mathcal{M}$ is cut out on
the~threefold $X$ by the~pencil
$$
\lambda x^{2}+\mu\Big(wx+q_{2}\big(x,y,z,t\big)\Big)=0,
$$
where $q_{i}$ and $p_{i}$ are  forms of degree $i$,  and
$(\lambda:\mu)\in\mathbb{P}^{1}$.
\end{itemize}
\end{theorem}

Let $P$ be an~arbitrary point of the~quartic hypersurface
$X\subset\mathbb{P}^{4}$.

\begin{definition}
\label{definition:threshold}  The mobility threshold of the~
threefold $X$ at the~point $P$ is the~number
$$
\iota\big(P\big)=\mathrm{sup}\left\{\lambda\in\mathbb{Q}\
\text{such that}\ \Big|n\Big(\pi^{*}\big(-K_{X}\big)-\lambda
E\Big)\Big|\ \text{has no fixed components for}\ n\gg 0\right\},
$$
where $\pi\colon Y\to X$ is the~ordinary blow up of $P$, and $E$
is the~exceptional divisor of $\pi$.
\end{definition}

Arguing as in the~proof of Theorem~\ref{theorem:main}, we obtain
the~following result.

\begin{theorem}
\label{theorem:threshold} The following conditions are equivalent:
\begin{itemize}
\item the~equality $\iota(P)=2$ holds,%
\item the~threefold $X$ can~be given by the~equation
$$
w^{3}x+w^{2}q_{2}\big(x,y,z,t\big)+wxp_{2}\big(x,y,z,t\big)+q_{4}\big(x,y,z,t\big)=0\subset\mathrm{Proj}\Big(\mathbb{C}[x,y,z,t,w]\Big)\cong\mathbb{P}^{4},%
$$
where $q_{i}$ and $p_{i}$ are forms of degree $i$ such that
$$
q_{2}\big(0,y,z,t\big)\ne 0,
$$
and $P$ is given by the~equations $x=y=z=t=0$.
\end{itemize}
\end{theorem}

One can~easily check that $2\geqslant\iota(P)\geqslant 1$.
Similarly, one can~show that
\begin{itemize}
\item $\iota(P)=1$ $\iff$ the~hyperplane section of $X$ that is singular at $P$ is a~cone,%
\item $\iota(P)=3/2$ $\iff$ the~threefold $X$ contains no lines passing through $P$.%
\end{itemize}

The proof of Theorem~\ref{theorem:main} is completed on board of
IL-96-300\,\emph{Valery~Chkalov}\,while\,flying from Seoul to
Moscow.
We\,thank\,Aeroflot\,Russian\,Airlines\,for\,good\,working\,conditions.

\section{Important lemma}
\label{section:preliminaries}

 Let $S$ be a~surface, let $O$ be a~smooth point of $S$, let $R$ be an~effective Weil divisor on
the~surface $S$, and let $\mathcal{D}$ be a~linear system on the~
surface $S$ that has no fixed components.

\begin{lemma}
\label{lemma:cornerstone} Let $D_{1}$ and $D_{2}$ be general
curves in $\mathcal{D}$. Then
$$
\mathrm{mult}_{O}\Big(D_{1}\cdot R\Big)=\mathrm{mult}_{O}\Big(D_{2}\cdot R\Big)\leqslant\mathrm{mult}_{O}\big(R\big)\mathrm{mult}_{O}\Big(D_{1}\cdot D_{1}\Big).%
$$
\end{lemma}

\begin{proof}
Put $S_{0}=S$ and $O_{0}=O$. Let us consider the~sequence of blow
ups
$$
\xymatrix{
&S_{n}\ar@{->}[rr]^{\pi_{n}}&&S_{n-1}\ar@{->}[rr]^{\pi_{n-1}}&&\cdots\ar@{->}[rr]^{\pi_{2}}&&S_{1}\ar@{->}[rr]^{\pi_{1}}&&S_{0}&}
$$
such that $\pi_{1}$ is a~blow up of the~point $O_{0}$, and
$\pi_{i}$ is a~blow up of the~point $O_{i-1}$ that is contained in
the~curve $E_{i-1}$, where $E_{i-1}$ is the~exceptional curve of
$\pi_{i-1}$, and $i=2,\ldots,n$.

Let $D^{i}_{j}$ be the~proper transform of $D_{j}$ on $S_{i}$ for
$i=0,\ldots,n$ and $j=1,2$. Then
$$
D^{i}_{1}\equiv D^{i}_{2}\equiv\pi_{i}^{*}\Big(D^{i-1}_{1}\Big)-\mathrm{mult}_{O_{i-1}}\Big(D^{i-1}_{1}\Big)E_{i}\equiv\pi_{i}^{*}\Big(D^{i-1}_{2}\Big)-\mathrm{mult}_{O_{i-1}}\Big(D^{i-1}_{2}\Big)E_{i}%
$$
for $i=1,\ldots,n$. Put
$d_{i}=\mathrm{mult}_{O_{i-1}}(D^{i-1}_{1})=\mathrm{mult}_{O_{i-1}}(D^{i-1}_{2})$
for $i=1,\ldots,n$.

Let $R^{i}$ be the~proper transform of $R$ on the~surface $S_{i}$
for $i=0,\ldots,n$. Then
$$
R^{i}\equiv\pi_{i}^{*}\Big(R^{i-1}\Big)-\mathrm{mult}_{O_{i-1}}\Big(R^{i-1}\Big)E_{i}
$$
for $i=1,\ldots,n$. Put $r_{i}=\mathrm{mult}_{O_{i-1}}(R^{i-1})$
for $i=1,\ldots,n$. Then $r_{1}=\mathrm{mult}_{O}(R)$.

We may chose the~blow ups $\pi_{1},\ldots,\pi_{n}$ in a~way such
that $D^{n}_{1}\cap D^{n}_{2}$ is empty in the~neighborhood of
the~exceptional locus of
$\pi_{1}\circ\pi_{2}\circ\cdots\circ\pi_{n}$. Then
$$
\mathrm{mult}_{O}\Big(D_{1}\cdot D_{2}\Big)=\sum_{i=1}^{n}d_{i}^{2}.%
$$

We may chose the~blow ups $\pi_{1},\ldots,\pi_{n}$ in a~way such
that $D^{n}_{1}\cap R^{n}$ and $D^{n}_{2}\cap R^{n}$ are empty in
the~neighborhood of the~exceptional locus of
$\pi_{1}\circ\pi_{2}\circ\cdots\circ\pi_{n}$. Then
$$
\mathrm{mult}_{O}\Big(D_{1}\cdot R\Big)=\mathrm{mult}_{O}\Big(D_{2}\cdot R\Big)=\sum_{i=1}^{n}d_{i}r_{i},%
$$
where some numbers among $r_{1},\ldots,r_{n}$ may be zero. Then
$$
\mathrm{mult}_{O}\Big(D_{1}\cdot R\Big)=\mathrm{mult}_{O}\Big(D_{2}\cdot R\Big)=\sum_{i=1}^{n}d_{i}r_{i}\leqslant\sum_{i=1}^{n}d_{i}r_{1}\leqslant\sum_{i=1}^{n}d^{2}_{i}r_{1}=\mathrm{mult}_{O}\big(R\big)\mathrm{mult}_{O}\Big(D_{1}\cdot D_{2}\Big),%
$$
because $d_{i}\leqslant d_{i}^{2}$ and $r_{i}\leqslant
r_{1}=\mathrm{mult}_{O}(R)$ for every $i=1,\ldots,n$.
\end{proof}

The assertion of Lemma~\ref{lemma:cornerstone} is a~cornerstone of
the~proof of Theorem~\ref{theorem:main}.

\section{Curves}
\label{section:curves}

Let $X$ be a~smooth quartic threefold in $\mathbb{P}^{4}$, let
$\mathcal{M}$ be a~Halphen pencil on $X$. Then
$$
\mathcal{M}\sim -nK_{X},
$$
since $\mathrm{Pic}(X)=\mathbb{Z}K_{X}$. Put $\mu=1/n$. Then
\begin{itemize}
\item the~log pair $(X, \mu\mathcal{M})$ is canonical by \cite[Theorem~A]{Isk01},%
\item  the~log pair  $(X, \mu\mathcal{M})$ is not terminal by \cite[Theorem~2.1]{ChPa05h}.%
\end{itemize}

Let $\mathbb{CS}(X, \mu\mathcal{M})$ be the~set of non-terminal
centers of  $(X, \mu\mathcal{M})$ (see~\cite{ChPa05h}).~Then
$$
\mathbb{CS}\Big(X, \mu\mathcal{M}\Big)\ne\varnothing,
$$
because $(X, \mu\mathcal{M})$ is not terminal. Let $M_{1}$ and
$M_{2}$ be two general surfaces in $\mathcal{M}$.

\begin{lemma}
\label{lemma:4n-square} Suppose that $\mathbb{CS}(X,
\mu\mathcal{M})$ contains a~point $P \in X$. Then
$$
\mathrm{mult}_{P}\big(M\big)=n\mathrm{mult}_{P}\big(T\big)=2n,
$$
where $M$ is any surface in $\mathcal{M}$, and $T$ is the~surface
in $|-K_{X}|$ that is singular at $P$.
\end{lemma}

\begin{proof}
It follows from \cite[Proposition~1]{Pu98a} that the~inequality
$$
\mathrm{mult}_{P}\Big(M_{1}\cdot M_{2}\Big)\geqslant 4n^{2}
$$
holds. Let $H$ be a~general surface in $|-K_{X}|$ such that $P\in
H$. Then
$$
4n^{2}=H\cdot M_{1}\cdot M_{2}\geqslant\mathrm{mult}_{P}\Big(M_{1}\cdot M_{2}\Big)\geqslant 4n^{2},%
$$
which gives $(M_{1}\cdot M_{2})_{P}=4n^{2}$. Arguing as in
the~proof of \cite[Proposition~1]{Pu98a}, we see that
$$
\mathrm{mult}_{P}\big(M_{1}\big)=\mathrm{mult}_{P}\big(M_{2}\big)=2n,
$$
because $(M_{1}\cdot M_{2})_{P}=4n^{2}$. Similarly, we see that
$$
4n=H\cdot T\cdot M_{1}\geqslant\mathrm{mult}_{P}\big(T\big)\mathrm{mult}_{P}\big(M_{1}\big)=2n\mathrm{mult}_{P}\big(T\big)\geqslant4n,%
$$
which implies that $\mathrm{mult}_{P}(T)=2$. Finally, we also have
$$
4n^{2}=H\cdot M\cdot M_{1}\geqslant\mathrm{mult}_{P}\big(M\big)\mathrm{mult}_{P}\big(M_{1}\big)=2n\mathrm{mult}_{P}\big(M\big)\geqslant4n^{2},%
$$
where $M$ is any surface in $\mathcal{M}$, which completes the~
proof.
\end{proof}

\begin{lemma}
\label{lemma:point-CS-lines} Suppose that $\mathbb{CS}(X,
\mu\mathcal{M})$ contains a~point $P\in X$.  Then
$$
M_{1}\cap M_{2}=\bigcup_{i=1}^{r}L_{i},
$$
where $L_{1},\ldots,L_{r}$ are lines on the~threefold $X$ that
pass through the~point $P$.
\end{lemma}

\begin{proof}
Let $H$ be a~general surface in $|-K_{X}|$ such that $P\in H$.
Then
$$
4n^{2}=H\cdot M_{1}\cdot M_{2}=\mathrm{mult}_{P}\Big(M_{1}\cdot M_{2}\Big)=4n^{2}%
$$
by Lemma~\ref{lemma:4n-square}. Then $\mathrm{Supp}(M_{1}\cdot
M_{2})$ consists of lines on $X$ that pass through $P$.
\end{proof}

\begin{lemma}
\label{lemma:lines-inequalities} Suppose that $\mathbb{CS}(X,
\mu\mathcal{M})$ contains a~point $P \in X$. Then
$$
n\big\slash 3\leqslant\mathrm{mult}_{L}\big(\mathcal{M}\big)\leqslant n\big\slash 2%
$$
for every line $L \subset X$ that passes through the~point $P$.
\end{lemma}

\begin{proof}
Let $D$ be a~general hyperplane section of $X$ through $L$. Then
we have
$$
M\Big\vert_{D}=\mathrm{mult}_{L}\big(\mathcal{M}\big)L+\Delta,
$$
where $M$ is a~general surface in $\mathcal{M}$ and $\Delta$ is
an~ effective divisor such that
$$
\mathrm{mult}_{P}\big(\Delta\big)\geqslant 2n-\mathrm{mult}_{L}\big(\mathcal{M}\big).%
$$

On the~surface $D$ we have $L\cdot L=-2$. Then
$$
n=\Big(\mathrm{mult}_{L}\big(\mathcal{M}\big)L + \Delta\Big) \cdot L= -2\mathrm{mult}_{L}\big(\mathcal{M}\big)+\Delta\cdot L%
$$
on the~surface $D$. But $\Delta~\cdot
L\geqslant\mathrm{mult}_{P}(\Delta)\geqslant
2n-\mathrm{mult}_{L}(\mathcal{M})$. Thus, we get
$$
n\geqslant -2\mathrm{mult}_{L}\big(\mathcal{M}\big) + \mathrm{mult}_{P}\big(\Delta\big)\geqslant 2n-3\mathrm{mult}_{L}\big(\mathcal{M}\big),%
$$
which implies that
$n/3\leqslant\mathrm{mult}_{L}\big(\mathcal{M}\big)$.

Let $T$ be the~surface in $|-K_{X}|$ that is singular at $P$. Then
$T\cdot D$ is reduced and
$$
T\cdot D = L+Z,
$$
where $Z$ is an~irreducible plane cubic curve such that $P\in Z$.
Then
$$
3n=\Big(\mathrm{mult}_{L}\big(\mathcal{M}\big)L + \Delta\Big)\cdot Z=3\mathrm{mult}_{L}\big(\mathcal{M}\big)+\Delta~\cdot Z%
$$
on the~surface $D$. The~set $\Delta\cap Z$ is finite by
Lemma~\ref{lemma:point-CS-lines}. In particular, we have
$$
\Delta~\cdot Z\geqslant\mathrm{mult}_{P}\big(\Delta\big)\geqslant 2n-\mathrm{mult}_{L}\big(\mathcal{M}\big),%
$$
because $\mathrm{Supp}(\Delta)$ does not contain the~curve $Z$.
Thus, we get
$$
3n\geqslant3\mathrm{mult}_{L}\big(\mathcal{M}\big)+\mathrm{mult}_{P}\big(\Delta\big)\geqslant
2n+2\mathrm{mult}_{L}\big(\mathcal{M}\big),
$$
which implies that
$\mathrm{mult}_{L}\big(\mathcal{M}\big)\leqslant n/2$.
\end{proof}

In the~rest of this section we prove the~following result.

\begin{proposition}
\label{proposition:curves} Suppose that $\mathbb{CS}(X,
\mu\mathcal{M})$ contains a~curve. Then $n=1$.
\end{proposition}

Suppose that $\mathbb{CS}(X, \mu\mathcal{M})$ contains a~curve
$Z$. Then it follows
Lemmas~\ref{lemma:point-CS-lines}~and~\ref{lemma:lines-inequalities}
that the~set $\mathbb{CS}(X, \mu\mathcal{M})$
does~not~contain~points of the~threefold $X$ and
\begin{equation}
\label{mult-of-M} \mathrm{mult}_{Z}\big(\mathcal{M}\big)=n,
\end{equation}
because $(X, \mu\mathcal{M})$ is canonical but not terminal. Then
$\mathrm{deg}(Z)\leqslant 4$ by \cite[Lemma~2.1]{ChPa05h}.

\begin{lemma}
\label{lemma:curve-is-line} Suppose that $\mathrm{deg}(Z)=1$. Then
$n = 1$.
\end{lemma}

\begin{proof}
Let $\pi: V \to X$ be the~blow up of $X$ along the~line $Z$.
Let~$\mathcal{B}$~be~the~proper~transform of the~pencil
$\mathcal{M}$ on the~threefold $V$, and let  $B$ be a general
surface in $\mathcal{B}$. Then
\begin{equation}
\label{mult-of-B} B \sim -nK_{V}
\end{equation}
by \eqref{mult-of-M}. There is a~commutative diagram
$$
\xymatrix{
&&V\ar@{->}[dl]_{\pi}\ar@{->}[dr]^{\eta}&&&\\%
&X\ar@{-->}[rr]_{\psi}&&\mathbb{P}^{2},&}
$$
where $\psi$ is the~projection from the~line $Z$ and $\eta$ is
the~ morphism induced by the~linear system $|-K_{V}|$. Thus, it
follows from \eqref{mult-of-B} that $\mathcal{B}$ is the~pull-back
of a pencil $\mathcal{P}$ on $\mathbb{P}^2$~by~$\eta$.

We see that the~base locus of $\mathcal{B}$ is contained in
the~union of fibers of $\eta$.

The~set $\mathbb{CS}(V, \mu\mathcal{B})$ is not empty by
\cite[Theorem~2.1]{ChPa05h}. It easily follows from
\eqref{mult-of-M} that the~set  $\mathbb{CS}(V, \mu\mathcal{B})$
does not contain points because $\mathbb{CS}(X, \mu\mathcal{M})$
contains no points.

We see that there is an~irreducible curve $L\subset V$ such that
$$
\mathrm{mult}_{L}\big(\mathcal{B}\big)=n
$$
and $\eta(L)$ is a~point $Q\in\mathbb{P}^{2}$. Let $C$ be a
general curve in $\mathcal{P}$. Then $\mathrm{mult}_{Q}(C)=n$. But
$$
C\sim \mathcal{O}_{\mathbb{P}^{2}}\big(n\big)
$$
by \eqref{mult-of-B}. Thus, we see that $n=1$, because general
surface in $\mathcal{M}$ is irreducible.
\end{proof}

Thus, we may assume that the~set $\mathbb{CS}(X, \mu\mathcal{M})$
does not contain lines.

\begin{lemma}
\label{lemma:plane-curve} The curve $Z\subset\mathbb{P}^{4}$ is
contained in a~plane.
\end{lemma}

\begin{proof}
Suppose that $Z$ is not contained in any plane in
$\mathbb{P}^{4}$. Let us show that this assumption leads to a
contradiction. Since $\mathrm{deg}(Z)\leqslant 4$, we have
$$
\mathrm{deg}\big(Z\big)\in\big\{3,4\big\},
$$
and $Z$ is smooth if $\mathrm{deg}(Z)=3$. If $\mathrm{deg}(Z)=4$,
then $Z$ may have at most one~double~point.

Suppose that $Z$ is smooth. Let $\alpha\colon U\to X$ be the~blow
up at $Z$, and let $F$ be the~exceptional divisor of the~morphism
$\alpha$. Then the~base locus of the~linear system
$$
\Big|\alpha^{*}\Big(-\mathrm{deg}\big(Z\big)K_{X}\Big)-F\Big|
$$
does not contain any curve. Let $D_{1}$ and $D_{2}$ be the~proper
transforms on $U$ of two~sufficiently  general surfaces in
the~linear system $\mathcal{M}$. Then it follows from
\eqref{mult-of-M} that
$$
\Bigg(\alpha^{*}\Big(-\mathrm{deg}\big(Z\big)K_{X}\Big)-F\Bigg)\cdot D_{1}\cdot D_{2}=n^{2}\Bigg(\alpha^{*}\Big(-\mathrm{deg}\big(Z\big)K_{X}\Big)-F\Bigg)\cdot\Bigg(\alpha^{*}\Big(-K_{X}\Big) - F\Bigg)^2 \geqslant 0,%
$$
because the~cycle $D_{1}\cdot D_{2}$ is effective. On the~other
hand, we have
$$
\Big(\alpha^{*}\Big(-\mathrm{deg}\big(Z\big)K_{X}\Big)-F\Big)\cdot\Big(\alpha^{*}\Big(-K_{X}\Big)-F\Big)^2=\Bigg(3\mathrm{deg}\big(Z\big)-\Big(\mathrm{deg}\big(Z\big)\Big)^{2}-2\Bigg)<0,%
$$
which is a contradiction. Thus, the~curve $Z$ is not smooth.

Thus, we see that $Z$ is a~quartic curve with a~double point $O$.

Let $\beta\colon W\to X$ be the~composition of the~blow up of
the~point $O$ with the~blow up of the~proper transform of
the~curve $Z$. Let $G$ and $E$ be the~exceptional surfaces of
the~morphism $\beta$ such that $\beta(E)=Z$ and $\beta(G)=O$. Then
the~base locus of the~linear~system
$$
\Big|\beta^{*}\big(-4K_{X}\big)-E-2G\Big|
$$
does not contain any curve. Let $R_{1}$ and $R_{2}$ be the~proper
transforms on $W$ of two~sufficiently general surfaces in
$\mathcal{M}$. Put $m=\mathrm{mult}_{O}(\mathcal{M})$. Then it
follows from \eqref{mult-of-M} that
$$
\Bigg(\beta^{*}\Big(-4K_{X}\Big)-E-2G\Bigg)\cdot R_{1}\cdot R_{2}=\Bigg(\beta^{*}\Big(-4K_{X}\Big)-E-2G\Bigg)\cdot\Bigg(\beta^{*}\Big(-nK_{X}\Big)-nE-mG\Bigg)^2\geqslant 0,%
$$
and $m<2n$, because the~set $\mathbb{CS}(X, \mu\mathcal{M})$ does
not contain points. Then
$$
\Bigg(\beta^{*}\Big(-4K_{X}\Big)-E-2G\Bigg)\cdot\Bigg(\beta^{*}\Big(-nK_{X}\Big)-nE-mG\Bigg)^2=-8n^2+6mn-m^2<0,
$$
which is a contradiction.
\end{proof}

If $\mathrm{deg}(Z)=4$, then $n=1$ by
Lemma~\ref{lemma:plane-curve} and \cite[Theorem~2.2]{ChPa05h}.

\begin{lemma}
\label{lemma:curve-is-cubic} Suppose that $\mathrm{deg}(Z)=3$.
Then $n=1$.
\end{lemma}

\begin{proof}
Let $\mathcal{P}$ be the~pencil in $|-K_{X}|$ that contains all
hyperplane sections of  $X$ that pass through the~curve $Z$. Then
the~base locus of $\mathcal{P}$ consists of the~curve $Z$ and a
line $L \subset X$.

Let $D$ be a~sufficiently general surface in the~pencil
$\mathcal{P}$, and let $M$ be a~sufficiently general surface in
the~pencil $\mathcal{M}$. Then $D$ is a~smooth surface, and
\begin{equation}
\label{on-D}
M\Big\vert_{D}=nZ+\mathrm{mult}_{L}\big(\mathcal{M}\big)L+B\equiv nZ+nL,%
\end{equation}
where $B$ is a curve whose support does not contain neither $Z$
nor $L$.

On the~surface $D$, we have $Z\cdot L=3$ and $L\cdot L=-2$.
Intersecting \eqref{on-D} with $L$, we get
$$
n=\big(nZ+nL\big)\cdot L=3n-2\mathrm{mult}_{L}\big(\mathcal{M}\big) + B \cdot L\geqslant 3n -2\mathrm{mult}_{L}\big(\mathcal{M}\big),%
$$
which easily implies that $\mathrm{mult}_{L}(\mathcal{M})
\geqslant n$. But the~inequality $\mathrm{mult}_{L}(\mathcal{M})
\geqslant n$ is impossible, because we assumed that
$\mathbb{CS}(X, \mu\mathcal{M})$ contains no lines.
\end{proof}

\begin{lemma}
\label{lemma:curve-is-conic} Suppose that $\mathrm{deg}(Z)=2$.
Then $n=1$.
\end{lemma}

\begin{proof}
Let $\alpha\colon U\to X$ be the~blow up of the~curve $Z$. Then
$|-K_{U}|$ is a pencil, whose base locus consists of a smooth
irreducible curve $L\subset U$.

Let $D$ be a~general surface in $|-K_{U}|$. Then $D$ is a~smooth
surface.

Let $\mathcal{B}$ be the~proper transform of the~pencil
$\mathcal{M}$ on the~threefold $U$. Then
$$
-nK_{U}\Big\vert_{D}\equiv B\Big\vert_{D}\equiv nL,
$$
where $B$ is a~general surface in $\mathcal{B}$. But $L^{2} = -2$
on the~surface $D$. Then
$$
L\in\mathbb{CS}\Big(U, \mu\mathcal{B}\Big)
$$
which implies that $\mathcal{B}=\big|-K_{U}\big|$ by
\cite[Theorem~2.2]{ChPa05h}. Then $n=1$.
\end{proof}

The assertion of Proposition~\ref{proposition:curves}~is~proved.

\section{Points}
\label{section:general-quartic}

Let $X$ be a~smooth quartic threefold in $\mathbb{P}^{4}$, let
$\mathcal{M}$ be a~Halphen pencil on $X$. Then
$$
\mathcal{M}\sim -nK_{X},
$$
since $\mathrm{Pic}(X)=\mathbb{Z}K_{X}$. Put $\mu=1/n$. Then
\begin{itemize}
\item the~log pair $(X, \mu\mathcal{M})$ is canonical by \cite[Theorem~A]{Isk01},%
\item  the~log pair  $(X, \mu\mathcal{M})$ is not terminal by \cite[Theorem~2.1]{ChPa05h}.%
\end{itemize}

\begin{remark}
\label{remark:n-not-1} To prove Theorem~\ref{theorem:main}, it is
enough to show that $X$ can~be given by
$$
w^{3}x+w^{2}q_{2}\big(x,y,z,t\big)+wxp_{2}\big(x,y,z,t\big)+q_{4}\big(x,y,z,t\big)=0\subset\mathrm{Proj}\Big(\mathbb{C}[x,y,z,t,w]\Big)\cong\mathbb{P}^{4},%
$$
where $q_{i}$ and $p_{i}$ are homogeneous polynomials of degree
$i\geqslant 2$ such that $q_{2}(0,y,z,t)\ne 0$.
\end{remark}

Let $\mathbb{CS}(X, \mu\mathcal{M})$ be the~set of non-terminal
centers of  $(X, \mu\mathcal{M})$ (see~\cite{ChPa05h}).~Then
$$
\mathbb{CS}\Big(X, \mu\mathcal{M}\Big)\ne\varnothing,
$$
because $(X, \mu\mathcal{M})$ is not terminal. Suppose that $n\ne
1$. There is a~point $P\in X$~such~that
$$
P\in \mathbb{CS}\Big(X,\mu\mathcal{M}\Big)
$$
by Proposition~\ref{proposition:curves}. It follows from
Lemmas~\ref{lemma:4n-square}, \ref{lemma:point-CS-lines} and
\ref{lemma:lines-inequalities} that
\begin{itemize}
\item there are finitely many distinct lines $L_{1},\ldots,L_{r}\subset X$ containing $P\in X$,%
\item the~equality  $\mathrm{mult}_{P}(M)=2n$ holds, and
$$
n/3\leqslant\mathrm{mult}_{L_{i}}\big(M\big)\leqslant n/2,
$$
where $M$ is a~general surface in the~pencil $\mathcal{M}$,%
\item the~equality $\mathrm{mult}_{P}(T)=2$ holds, where $T\in|-K_{X}|$ such that $\mathrm{mult}_{P}(T)\geqslant 2$,%
\item the~base locus of the~pencil $\mathcal{M}$ consists of
the~lines $L_{1},\ldots,L_{r}$, and
$$
\mathrm{mult}_{P}\Big(M_{1}\cdot M_{2}\Big)=4n^{2},
$$
where $M_{1}$ and $M_{2}$ are sufficiently general surfaces in
$\mathcal{M}$.
\end{itemize}

\begin{lemma}
\label{lemma:point-is-unique} The equality
$\mathbb{CS}(X,\mu\mathcal{M})=\{P\}$ holds.
\end{lemma}

\begin{proof}
The~set $\mathbb{CS}(X, \mu\mathcal{M})$ does not contain curves
by Proposition~\ref{proposition:curves}.

Suppose that $\mathbb{CS}(X, \mu\mathcal{M})$ contains a~point
$Q\in X$ such that $Q\ne P$. Then $r=1$.

Let $D$ be a~general hyperplane section of $X$ that passes through
$L_{1}$. Then
$$
M\Big\vert_{D}=\mathrm{mult}_{L_{1}}\big(\mathcal{M}\big)L_{1}+\Delta,
$$
where $M$ is a~general surface in $\mathcal{M}$ and $\Delta$ is
an~ effective divisor such that
$$
\mathrm{mult}_{P}\big(\Delta\big)\geqslant 2n-\mathrm{mult}_{L_{1}}\big(\mathcal{M}\big)\leqslant\mathrm{mult}_{Q}\big(\Delta\big).%
$$

On the~surface $D$, we have $L_{1}^2=-2$. Then
$$
n=\Big(\mathrm{mult}_{L_{1}}\big(\mathcal{M}\big)L_{1} +
\Delta\Big) \cdot
L_{1}=-2\mathrm{mult}_{L_{1}}\big(\mathcal{M}\big) + \Delta\cdot
L\geqslant
-2\mathrm{mult}_{L_{1}}\big(\mathcal{M}\big)+2\Big(2n-\mathrm{mult}_{L_{1}}\big(\mathcal{M}\big)\Big),
$$
which gives $\mathrm{mult}_{L_{1}}(\mathcal{M})\geqslant 3n/4$.
But $\mathrm{mult}_{L_{1}}(\mathcal{M})\leqslant n/2$ by
Lemma~\ref{lemma:lines-inequalities}.
\end{proof}

The quartic threefold $X$ can~be given by an~equation
$$
w^{3}x+w^{2}q_{2}\big(x,y,z,t\big)+wq_{3}\big(x,y,z,t\big)+q_{4}\big(x,y,z,t\big)=0\subset\mathrm{Proj}\Big(\mathbb{C}[x,y,z,t,w]\Big)\cong\mathbb{P}^{4},%
$$
where $q_{i}$ is a~homogeneous polynomial of degree $i\geqslant
2$.

\begin{remark}
\label{remark:lines-algebraic} The lines
$L_{1},\ldots,L_{r}\subset\mathbb{P}^{4}$ are given by the~
equations
$$
x=q_{2}\big(x,y,z,t\big)=q_{3}\big(x,y,z,t\big)=q_{4}\big(x,y,z,t\big)=0,
$$
the~surface  $T$ is cut out on $X$ by $x=0$, and
$\mathrm{mult}_{P}(T)=2\iff q_{2}(0,y,z,t)\ne 0$.
\end{remark}

Let $\pi\colon V\to X$ be the~blow up of the~point $P$, let $E$ be
the~$\pi$-exceptional divisor. Then
$$
\mathcal{B}\equiv\pi^{*}\big(-nK_{X}\big)-2nE\equiv -nK_{V},
$$
where $\mathcal{B}$ is the~proper transform of the~pencil
$\mathcal{M}$ on the~threefold $V$.

\begin{remark}
\label{remark:no-new-base-curves} The pencil $\mathcal{B}$ has no
base curves in $E$, because
$$
\mathrm{mult}_{P}\big(M_{1}\cdot M_{2}\big)=\mathrm{mult}_{P}\big(M_{1}\big)\mathrm{mult}_{P}\big(M_{2}\big).%
$$
\end{remark}

Let $\bar{L}_{i}$ be the~proper transform of the~line $L_{i}$ on
the~threefold $V$ for $i=1,\ldots,r$. Then
$$
B_{1}\cdot B_{2}=\sum_{i=1}^{r}\mathrm{mult}_{\bar{L}_{i}}\Big(B_{1}\cdot B_{2}\Big)\bar{L}_{i},%
$$
where $B_{1}$ and $B_{2}$ are proper transforms of $M_{1}$ and
$M_{2}$ on the~threefold $V$, respectively.

\begin{lemma}
\label{lemma:curves-2mult-deg} Let $Z$ be an~irreducible curve on
$X$ such that $Z\not\in\{L_{1},\ldots,Z_{r}\}$. Then
$$
\mathrm{deg}\big(Z\big)\geqslant 2\mathrm{mult}_{P}\big(Z\big),
$$
and the~equality $\mathrm{deg}(Z)=2\mathrm{mult}_{P}(Z)$ implies
that
$$
\bar{Z}\cap\left(\bigcup_{i=1}^{r}\bar{L}_{i}\right)=\varnothing,
$$
where $\bar{Z}$ is a~proper transform of the~curve $Z$ on the~
threefold $V$.
\end{lemma}

\begin{proof}
The curve $\bar{Z}$ is not contained in the~base locus of
the~pencil $\mathcal{B}$. Then
$$
0\leqslant B_{i}\cdot \bar{Z}\leqslant
n\Big(\mathrm{deg}\big(Z\big)-2\mathrm{mult}_{P}\big(Z\big)\Big),
$$
which implies the~required assertions.
\end{proof}

To conclude the~proof of Theorem~\ref{theorem:main}, it is enough
to show that
$$
q_{3}\big(x,y,z,t\big)=xp_{2}\big(x,y,z,t\big)+q_{2}\big(x,y,z,t\big)p_{1}\big(x,y,z,t\big),
$$
where $p_{1}$ and $p_{2}$ are some homogeneous polynomials of
degree $1$ and $2$, respectively.

\section{Good points}
\label{section:good-points}

Let us use the~assumptions and notation of
Section~\ref{section:general-quartic}. Suppose that the~conic
$$
q_{2}(0,y,z,t)=0\subset\mathrm{Proj}\Big(\mathbb{C}[y,z,t]\Big)\cong\mathbb{P}^{2}
$$
is reduced and irreducible. In this section we prove the~following
result.

\begin{proposition}
\label{proposition:good-points} The polynomial $q_{3}(0,y,z,t)$ is
divisible by $q_{2}(0,y,z,t)$.
\end{proposition}

Let us prove Proposition~\ref{proposition:good-points}. Suppose
that $q_{3}(0,y,z,t)$ is not divisible by $q_{2}(0,y,z,t)$.

Let $\mathcal{R}$ be the~linear system on the~threefold $X$ that
is cut out by quadrics
$$
xh_{1}(x,y,z,t)+\lambda \Big(wx+q_{2}(x,y,z,t)\Big)=0,
$$
where $h_{1}$ is an~arbitrary linear form and
$\lambda\in\mathbb{C}$. Then $\mathcal{R}$ does not have fixed
components.

\begin{lemma}
\label{lemma:ODP-2H-3E} Let $R_{1}$ and $R_{2}$ be general
surfaces in the~linear system $\mathcal{R}$. Then
$$
\sum_{i=1}^{r}\mathrm{mult}_{L_{i}}\Big(R_{1}\cdot R_{2}\Big)\leqslant 6.%
$$
\end{lemma}

\begin{proof}
We may assume that $R_{1}$ is cut out by the~equation
$$
wx+q_{2}(x,y,z,t)=0,
$$
and $R_{2}$ is cut out by $xh_{1}(x,y,z,t)=0$, where $h_{1}$ is
sufficiently general. Then
$$
\mathrm{mult}_{L_{i}}\Big(R_{1}\cdot
R_{2}\Big)=\mathrm{mult}_{L_{i}}\Big(R_{1}\cdot T\Big).
$$

Put $m_{i}=\mathrm{mult}_{L_{i}}(R_{1}\cdot T)$. Then
$$
R_{1}\cdot T=\sum_{i=1}^{r}m_{i}L_{i}+\Delta,
$$
where $m_{i}\in\mathbb{N}$, and $\Delta$ is a~cycle, whose support
contains no lines passing through $P$.

Let $\bar{R}_{1}$ and $\bar{T}$ be the~proper transforms of
$R_{1}$ and $T$ on $V$, respectively. Then
$$
\bar{R}_{1}\cdot\bar{T}=\sum_{i=1}^{r}m_{i}\bar{L}_{i}+\Omega,
$$
where $\Omega$ is an~effective cycle, whose support contains no
lines passing through $P$.

The support of the~cycle $\Omega$ does not contain curves that are
contained in the~exceptional divisor $E$, because $q_{3}(0,y,z,t)$
is not divisible by $q_{2}(0,y,z,t)$ by our assumption. Then
$$
6=E\cdot\bar{R}_{1}\cdot\bar{T}=\sum_{i=1}^{r}m_{i}\big(E\cdot\bar{L}_{i}\big)+E\cdot\Omega\geqslant
\sum_{i=1}^{r}m_{i}\big(E\cdot\bar{L}_{i}\big)=\sum_{i=1}^{r}m_{i},
$$
which is exactly what we want.
\end{proof}

Let $M$ and $R$ be general surfaces in $\mathcal{M}$ and
$\mathcal{R}$, respectively. Put
$$
M\cdot R=\sum_{i=1}^{r}m_{i}L_{i}+\Delta,
$$
where $m_{i}\in\mathbb{N}$, and $\Delta$ is a~cycle, whose support
contains no lines  passing through $P$.

\begin{lemma}
\label{lemma:good-points-Delta-not-empty} The cycle $\Delta$ is
not trivial.
\end{lemma}

\begin{proof}
Suppose that $\Delta=0$. Then $\mathcal{M}=\mathcal{R}$ by
\cite[Theorem~2.2]{ChPa05h}. But $\mathcal{R}$ is not a~pencil.
\end{proof}

We have $\mathrm{deg}(\Delta)=8n-\sum_{i=1}^{r} m_{i}$. On the~
other hand, the~inequality
$$
\mathrm{mult}_{P}\big(\Delta\big)\geqslant 6n-\sum_{i=1}^{r} m_{i}
$$
holds, because $\mathrm{mult}_{P}(M)=2n$ and
$\mathrm{mult}_{P}(R)\geqslant 3$. It follows from
Lemma~\ref{lemma:curves-2mult-deg} that
$$
\mathrm{deg}\big(\Delta\big)=8n-\sum_{i=1}^{r} m_{i}\geqslant 2\mathrm{mult}_{P}\big(\Delta\big)\geqslant 2\left(6n-\sum_{i=1}^{r} m_{i}\right),%
$$
which implies that $\sum_{i=1}^{r} m_{i}\geqslant 4n$. But it
follows from Lemmas~\ref{lemma:cornerstone} and
\ref{lemma:lines-inequalities} that
$$
m_{i}\leqslant \mathrm{mult}_{L_{i}}\Big(R_{1}\cdot R_{2}\Big)\mathrm{mult}_{L_{i}}\big(M\big)\leqslant \mathrm{mult}_{L_{i}}\Big(R_{1}\cdot R_{2}\Big)n/2%
$$
for every $i=1,\ldots,r$, where $R_{1}$ and $R_{2}$ are general
surfaces in $\mathcal{R}$. Then
$$
\sum_{i=1}^{r}m_{i}\leqslant
\sum_{i=1}^{r}\mathrm{mult}_{L_{i}}\Big(R_{1}\cdot R_{2}\Big)n/2\leqslant 3n%
$$
by Lemma~\ref{lemma:ODP-2H-3E}, which is a~contradiction.

The assertion of Proposition~\ref{proposition:good-points} is
proved.

\section{Bad points}
\label{section:bad-points}

Let us use the~assumptions and notation of
Section~\ref{section:general-quartic}. Suppose that the~conic
$$
q_{2}(0,y,z,t)=0\subset\mathrm{Proj}\Big(\mathbb{C}[y,z,t]\Big)\cong\mathbb{P}^{2}
$$
is reduced and reducible. Therefore, we have
$$
q_{2}(x,y,z,t)=\big(\alpha_{1}y+\beta_{1}z+\gamma_{1}t\big)\big(\alpha_{2}y+\beta_{2}z+\gamma_{2}t\big)+xp_{1}\big(x,y,z,t\big)
$$
where $p_{1}(x,y,z,t)$ is a~linear form, and
$(\alpha_{1}:\beta_{1}:\gamma_{1})\in\mathbb{P}^{2}\ni
(\alpha_{2}:\beta_{2}:\gamma_{2})$.

\begin{proposition}
\label{proposition:bad-points} The polynomial $q_{3}(0,y,z,t)$ is
divisible by $q_{2}(0,y,z,t)$.
\end{proposition}

Suppose that $q_{3}(0,y,z,t)$ is not divisible by
$q_{2}(0,y,z,t)$. Then without loss of generality, we may assume
that $q_{3}(0,y,z,t)$ is not divisible by
$\alpha_{1}y+\beta_{1}z+\gamma_{1}t$.

Let $Z$ be the~curve in $X$ that is cut out by the~equations
$$
x=\alpha_{1}y+\beta_{1}z+\gamma_{1}t=0.
$$

\begin{remark}
\label{remark:bad-points-Z-curve-not-reduced} The equality
$\mathrm{mult}_{P}(Z)=3$ holds, but $Z$ is not necessary reduced.
\end{remark}

Hence, it follows from Lemma~\ref{lemma:curves-2mult-deg} that
$\mathrm{Supp}(Z)$ contains a~line among $L_{1},\ldots,L_{r}$.

\begin{lemma}
\label{lemma:bad-points-no-conics} The support of the~curve $Z$
does not contain an~irreducible conic.
\end{lemma}

\begin{proof}
Suppose that $\mathrm{Supp}(Z)$ contains an~irreducible conic $C$.
Then
$$
Z=C+L_{i}+L_{j}
$$
for some $i\in\{1,\ldots,r\}\ni j$. Then $i=j$, because otherwise
the~set
$$
\Big(C\cap L_{i}\Big)\bigcup\Big(C\cap L_{j}\Big)
$$
contains a~point that is different from $P$, which is impossible
by Lemma~\ref{lemma:curves-2mult-deg}. We see that
$$
Z=C+2L_{i},
$$
and it follows from Lemma~\ref{lemma:curves-2mult-deg} that $C\cap
L_{i}=P$. Then $C$ is tangent to $L_{i}$ at the~point $P$

Let $\bar{C}$ be a~proper transform of the~curve $C$ on the~
threefold $V$. Then
$$
\bar{C}\cap\bar{L}_{i}\ne\varnothing,
$$
which is impossible by Lemma~\ref{lemma:curves-2mult-deg}. The
assertion is proved.
\end{proof}

\begin{lemma}
\label{lemma:bad-points-Z-is-4-lines} The support of the~curve $Z$
consists of lines.
\end{lemma}
\begin{proof}
Suppose that $\mathrm{Supp}(Z)$ does not consist of lines. It
follows from Lemma~\ref{lemma:bad-points-no-conics} that
$$
Z=L_{i}+C,
$$
where $C$ is an~irreducible cubic curve. But
$\mathrm{mult}_{P}(Z)=3$. Then
$$
\mathrm{mult}_{P}\big(C\big)=2,
$$
which is impossible by
Lemma~\ref{lemma:curves-2mult-deg}
\end{proof}

We may assume that there is a line $L\subset X$ such that
$P\not\in P$ and
$$
Z=a_{1}L_{1}+\cdots+a_{k}L_{k}+L,
$$
where $a_{1},a_{2},a_{3}\in\mathbb{N}$ such that $a_{1}\geqslant
a_{2}\geqslant a_{3}$ and $\sum_{i=1}^{k}a_{i}=3$.

\begin{remark}
\label{remark:bad-points-counting-lines} We have $L_{i}\ne L_{j}$
whenever $i\ne j$.
\end{remark}

Let $H$ be a~sufficiently general surface of $X$ that is cut out
by the~equation
$$
\lambda x+\mu\Big(\alpha_{1}y+\beta_{1}z+\gamma_{1}t\Big)=0,
$$
where $(\lambda:\mu)\in\mathbb{P}^{1}$. Then $H$ has at most
isolated singularities.

\begin{remark}
\label{remark:bad-points-normal-surface} The surface $H$ is smooth
at the~points $P$ and $L\cap L_{i}$, where $i=1,\ldots,k$.
\end{remark}

 Let $\bar{H}$
and $\bar{L}$ be the~proper transforms of $H$ and $L$ on the~
threefold $V$, respectively.

\begin{lemma}
\label{lemma:bad-points-3-lines} The inequality $k\ne 3$ holds.
\end{lemma}

\begin{proof}
Suppose that the~equality $k=3$ holds. Then $H$ is smooth. Put
$$
B\Big\vert_{\bar{H}}=m_{1}\bar{L}_{1}+m_{2}\bar{L}_{2}+m_{3}\bar{L}_{3}+\Omega,
$$
where $B$ is a~general surface in $\mathcal{B}$, and $\Omega$ is
an~effective divisor on $\bar{H}$ whose support does not contain
any of the~curves $\bar{L}_{1}$, $\bar{L}_{2}$ and $\bar{L}_{3}$.
Then
$$
\bar{L}\not\subseteq\mathrm{Supp}\big(\Omega\big)\not\supseteq\bar{H}\cap E,%
$$
because the~base locus of the~pencil $\mathcal{B}$ consists of
the~curves $\bar{L}_{1},\ldots,\bar{L}_{r}$. Then
$$
n=\bar{L}\cdot\Big(m_{1}\bar{L}_{1}+m_{2}\bar{L}_{2}+m_{3}\bar{L}_{3}+\Omega\Big)=\sum_{i=1}^{3} m_{i}+\bar{L}\cdot\Omega\geqslant \sum_{i=1}^{3} m_{i},%
$$
which implies that $\sum_{i=1}^{3} m_{i}\leqslant n$. On the~other
hand, we have
$$
-n=\bar{L}_{i}\cdot\Big(m_{1}\bar{L}_{1}+m_{2}\bar{L}_{2}+m_{3}\bar{L}_{3}+\Omega\Big)=-3m_{i}+L_{i}\cdot\Omega\geqslant -3m_{i},%
$$
which implies that $m_{i}\geqslant n/3$. Thus, we have
$m_{1}=m_{2}=m_{3}=n/3$ and
$$
\Omega\cdot \bar{L}=\Omega\cdot\bar{L}_{1}=\Omega\cdot\bar{L}_{2}=\Omega\cdot\bar{L}_{3}=0,%
$$
which implies that
$\mathrm{Supp}(\Omega)\cap\bar{L}_{1}=\mathrm{Supp}(\Omega)\cap\bar{L}_{2}=\mathrm{Supp}(\Omega)\cap\bar{L}_{3}=\varnothing$.

Let $B^{\prime}$ be another general surface in $\mathcal{B}$.
Arguing as above, we see that
$$
B^{\prime}\Big\vert_{\bar{H}}=\frac{n}{3}\Big(\bar{L}_{1}+\bar{L}_{2}+\bar{L}_{3}\Big)+\Omega^{\prime},%
$$
where $\Omega^{\prime}$ is an~effective divisor on the~surface
$\bar{H}$ such~that
$$
\mathrm{Supp}\big(\Omega^{\prime}\big)\cap\bar{L}_{1}=\mathrm{Supp}\big(\Omega^{\prime}\big)\cap\bar{L}_{2}=\mathrm{Supp}\big(\Omega^{\prime}\big)\cap\bar{L}_{3}=\varnothing.
$$

One can easily check that $\Omega\cdot\Omega^{\prime}=n^{2}\ne 0$.
Then
$$
\mathrm{Supp}\big(\Omega\big)\cap\mathrm{Supp}\big(\Omega^{\prime}\big)\ne\varnothing,
$$
because
$|\mathrm{Supp}(\Omega)\cap\mathrm{Supp}(\Omega^{\prime})|<+\infty$
due to generality of the~surfaces $B$ and $B^{\prime}$.

The base locus of the~pencil $\mathcal{B}$ consists of the~curves
$\bar{L}_{1},\ldots,\bar{L}_{r}$. Hence, we have
$$
\mathrm{Supp}\big(B\big)\cap\mathrm{Supp}\big(B^{\prime}\big)=\bigcup_{i=1}^{r}\bar{L}_{i},
$$
but $\bar{L}_{i}\cap \bar{H}=\varnothing$ whenever
$i\not\in\{1,2,3\}$. Hence, we have
$$
\bar{L}_{1}\cup\bar{L}_{2}\cup\bar{L}_{3}\cup\left(\mathrm{Supp}\big(\Omega\big)\cap\mathrm{Supp}\big(\Omega^{\prime}\big)\right)
=\mathrm{Supp}\big(B\big)\cap\mathrm{Supp}\big(B^{\prime}\big)\cap\bar{H}=\bar{L}_{1}\cup\bar{L}_{2}\cup\bar{L}_{3},
$$
which implies that
$\mathrm{Supp}(\Omega)\cap\mathrm{Supp}(\Omega^{\prime})\subset\bar{L}_{1}\cup\bar{L}_{2}\cup\bar{L}_{3}$.
In particular, we see that
$$
\mathrm{Supp}(\Omega)\cap\Big(\bar{L}_{1}\cup\bar{L}_{2}\cup\bar{L}_{3}\Big)\ne\varnothing,
$$
because
$\mathrm{Supp}(\Omega)\cap\mathrm{Supp}(\Omega^{\prime})\ne\varnothing$.
But $\mathrm{Supp}(\Omega)\cap\bar{L}_{i}=\varnothing$ for
$i=1,2,3$.
\end{proof}

\begin{lemma}
\label{lemma:bad-points-2-lines} The inequality $k\ne 2$ holds.
\end{lemma}

\begin{proof}
Suppose that the~equality $k=2$ holds. Then $Z=2L_{1}+L_{2}+L$.
Put
$$
B\Big\vert_{\bar{H}}=m_{1}\bar{L}_{1}+m_{2}\bar{L}_{2}+\Omega,
$$
where $B$ is a~general surface in $\mathcal{B}$, and $\Omega$ is
an~effective divisor on  $\bar{H}$ whose support does not contain
the~curves $\bar{L}_{1}$ and $\bar{L}_{2}$. Then
$\bar{L}\not\subseteq\mathrm{Supp}(\Omega)\not\supseteq\bar{H}\cap
E$ and
$$
n=\bar{L}\cdot\Big(m_{1}\bar{L}_{1}+m_{2}\bar{L}_{2}+\Omega\Big)=m_{1}+m_{2}+\bar{L}\cdot\Omega\geqslant m_{1}+m_{2},%
$$
which implies that $m_{1}+m_{2}\leqslant n$. On the~other hand, we
have
$$
\bar{T}\vert_{\bar{H}}=2\bar{L}_{1}+\bar{L}_{2}+\bar{L}+E\Big\vert_{\bar{H}}\equiv\left(\pi^{*}\Big(-K_{X}\Big)-2E\right)\Big\vert_{\bar{H}},
$$
where $\bar{T}$ is the~proper transform of the~surface $T$ on the~
threefold $V$. Then
$$
-1=\bar{L}_{1}\cdot\Big(2\bar{L}_{1}+\bar{L}_{2}+\bar{L}+E\Big\vert_{\bar{H}}\Big)=2\Big(\bar{L}_{1}\cdot\bar{L}_{1}\Big)+2,
$$
which implies that $\bar{L}_{1}\cdot\bar{L}_{1}=-3/2$ on the~
surface $\bar{H}$. Then
$$
-n=\bar{L}_{1}\cdot\Big(m_{1}\bar{L}_{1}+m_{2}\bar{L}_{2}+\Omega\Big)=-3m_{1}/2+L_{1}\cdot\Omega\geqslant -3m_{1}/2,%
$$
which gives $m_{1}\geqslant 2n/3$. Similarly, we see that
$\bar{L}_{2}\cdot\bar{L}_{2}=-3$ on the~surface $\bar{H}$. Then
$$
-n=\bar{L}_{2}\cdot\Big(m_{1}\bar{L}_{1}+m_{2}\bar{L}_{2}+\Omega\Big)=-3m_{2}+L_{2}\cdot\Omega\geqslant -3m_{2},%
$$
which implies that $m_{2}\leqslant n/3$. Thus, we have
$m_{1}=2m_{2}=2n/3$ and
$$
\Omega\cdot \bar{L}=\Omega\cdot\bar{L}_{1}=\Omega\cdot\bar{L}_{2}=0,%
$$
which implies that
$\mathrm{Supp}(\Omega)\cap\bar{L}_{1}=\mathrm{Supp}(\Omega)\cap\bar{L}_{2}=\varnothing$.

Let $B^{\prime}$ be another general surface in $\mathcal{B}$.
Arguing as above, we see that
$$
B^{\prime}\Big\vert_{\bar{H}}=\frac{2n}{3}\bar{L}_{1}+\frac{n}{3}\bar{L}_{2}+\Omega^{\prime},%
$$
where $\Omega^{\prime}$ is an~effective divisor on $\bar{H}$ whose
support does not contain $\bar{L}_{1}$ and $\bar{L}_{2}$ such that
$$
\mathrm{Supp}\big(\Omega^{\prime}\big)\cap\bar{L}_{1}=\mathrm{Supp}\big(\Omega^{\prime}\big)\cap\bar{L}_{2}=\varnothing,
$$
which implies that $\Omega\cdot\Omega^{\prime}=n^{2}$. In
particular, we see that
$$
\mathrm{Supp}\big(\Omega\big)\cap\mathrm{Supp}\big(\Omega^{\prime}\big)\ne\varnothing,
$$
and arguing as in the~proof of
Lemma~\ref{lemma:bad-points-3-lines} we obtain a~contradiction.
\end{proof}

It follows from Lemmas~\ref{lemma:bad-points-3-lines} and
\ref{lemma:bad-points-2-lines} that $Z=3L_{1}+L$. Put
$$
B\Big\vert_{\bar{H}}=m_{1}\bar{L}_{1}+\Omega,
$$
where $B$ is a~general surface in $\mathcal{B}$, and $\Omega$ is
a~curve such that $\bar{L}_{1}\not\subseteq\mathrm{Supp}(\Omega)$.
Then
$$
\bar{L}\not\subseteq\mathrm{Supp}\big(\Omega\big)\not\supseteq\bar{H}\cap E,%
$$
because the~base locus of $\mathcal{B}$ consists of the~curves
$\bar{L}_{1},\ldots,\bar{L}_{r}$. Then
$$
n=\bar{L}\cdot\Big(m_{1}\bar{L}_{1}+\Omega\Big)=m_{1}+\bar{L}\cdot\Omega\geqslant m_{1},%
$$
which implies that $m_{1}\leqslant n$. On the~other hand, we have
$$
\bar{T}\vert_{\bar{H}}=3\bar{L}_{1}+\bar{L}+E\Big\vert_{\bar{H}}\equiv\left(\pi^{*}\Big(-K_{X}\Big)-2E\right)\Big\vert_{\bar{H}},
$$
where $\bar{T}$ is the~proper transform of the~surface $T$ on the~
threefold $V$. Then
$$
-1=\bar{L}_{1}\cdot\Big(3\bar{L}_{1}+\bar{L}+E\Big\vert_{\bar{H}}\Big)=3\bar{L}_{1}\cdot\bar{L}_{1}+2,
$$
which implies that $\bar{L}_{1}\cdot\bar{L}_{1}=-1$ on the~surface
$\bar{H}$. Then
$$
-n=\bar{L}_{1}\cdot\Big(m_{1}\bar{L}_{1}+\Omega\Big)=-m_{1}+L_{1}\cdot\Omega\geqslant -m_{1},%
$$
which gives $m_{1}\geqslant n$. Thus, we have $m_{1}=n$ and
$\Omega\cdot \bar{L}=\Omega\cdot\bar{L}_{1}=0$. Then
$\mathrm{Supp}(\Omega)\cap\bar{L}_{1}=\varnothing$.

Let $B^{\prime}$ be another general surface in $\mathcal{B}$.
Arguing as above, we see that
$$
B^{\prime}\Big\vert_{\bar{H}}=n\bar{L}_{1}+\Omega^{\prime},%
$$
where $\Omega^{\prime}$ is an~effective divisor on $\bar{H}$ whose
support does not contain $\bar{L}_{1}$ such that
$$
\mathrm{Supp}\big(\Omega^{\prime}\big)\cap\bar{L}_{1}=\varnothing,
$$
which implies that $\Omega\cdot\Omega^{\prime}=n^{2}$. In
particular, we see that
$\mathrm{Supp}(\Omega)\cap\mathrm{Supp}(\Omega^{\prime})\ne\varnothing$.

The base locus of the~pencil $\mathcal{B}$ consists of the~curves
$\bar{L}_{1},\ldots,L_{r}$. Hence, we have
$$
\mathrm{Supp}\big(B\big)\cap\mathrm{Supp}\big(B^{\prime}\big)=\bigcup_{i=1}^{r}\bar{L}_{i},
$$
but $\bar{L}_{i}\cap \bar{H}=\varnothing$ whenever
$\bar{L}_{i}\ne\bar{L}_{1}$. Then
$\mathrm{Supp}(\Omega)\cap\bar{L}_{1}\ne\varnothing$, because
$$
\bar{L}_{1}\cup\Big(\mathrm{Supp}\big(\Omega\big)\cap\mathrm{Supp}\big(\Omega^{\prime}\big)\Big)
=\mathrm{Supp}\big(B\big)\cap\mathrm{Supp}\big(B^{\prime}\big)\cap\bar{H}=\bar{L}_{1},
$$
which is a~contradiction. The assertion of
Proposition~\ref{proposition:bad-points} is proved.

\section{Very bad points}
\label{section:very-bad-points}

Let us use the~assumptions and notation of
Section~\ref{section:general-quartic}. Suppose that $q_{2}=y^{2}$.

The proof of Proposition 6.1 implies that $q_{3}(0,y,z,t)$ is
divisible by $y$. Then
$$
q_{3}=yf_{2}\big(z,t\big)+xh_{2}\big(z,t\big)+x^{2}a_{1}\big(x,y,z,t\big)+xyb_{1}\big(x,y,z,t\big)+y^{2}c_{1}\big(y,z,t\big)%
$$
where $a_{1}$, $b_{1}$, $c_{1}$ are linear forms, $f_{2}$ and
$h_{2}$ is are homogeneous polynomials of degree two.

\begin{proposition}
\label{proposition:very-bad-points} The equality $f_{2}(z,t)=0$
holds.
\end{proposition}

Let us prove Proposition~\ref{proposition:very-bad-points} by
reductio ad absurdum. Suppose that $f_{2}(z,t)\ne 0$.

\begin{remark}
\label{remark:zt-or-z2} By choosing suitable coordinates, we may
assume that $f_{2}=zt$ or $f_{2}=z^{2}$.
\end{remark}

We must use smoothness of the~threefold $X$ by analyzing the~shape
of $q_{4}$.~We~have
$$
q_{4}=f_{4}\big(z,t\big)+xu_{3}\big(z,t\big)+yv_{3}\big(z,t\big)+x^{2}a_{2}(x,y,z,t)+xyb_{2}(x,y,z,t)+y^{2}c_{2}(y,z,t),
$$
where $a_{2}$, $b_{2}$, $c_{2}$ are homogeneous polynomials of
degree two, $u_{3}$ and $v_{3}$ are homogeneous polynomials of
degree three, and $f_{4}$ is a~homogeneous polynomial of degree
four.

\begin{lemma}
\label{lemma:very-bad-points-derivatives} Suppose that
$f_{2}(z,t)=zt$ and
$$
f_{4}\big(z,t\big)=t^{2}g_{2}\big(z,t\big)
$$
for some $g_{2}(z,t)\in\mathbb{C}[z,t]$. Then $v_{3}(z,0)\ne 0$.
\end{lemma}

\begin{proof}
Suppose that $v_{3}(z,0)=0$. The surface $T$ is given by
the~equation
$$
w^{2}y^{2}+yzt+y^{2}c_{1}\big(x,y,z,t\big)+t^{2}g_{2}\big(z,t\big)+yv_{3}\big(z,t\big)+y^{2}c_{2}\big(x,y,z,t\big)
=
0\subset\mathrm{Proj}\Big(\mathbb{C}[y,z,t,w]\Big)\cong\mathbb{P}^{3}
$$
because $T$ is cut out on $X$ by the~equation $x=0$. Then $T$ has
non-isolated singularity along the~line $x=y=t=0$, which is
impossible because $X$ is smooth.
\end{proof}

Arguing as in the~proof of
Lemma~\ref{lemma:very-bad-points-derivatives}, we obtain
the~following corollary.

\begin{corollary}
\label{corollary:very-bad-points-derivatives} Suppose that
$f_{2}(z,t)=zt$ and
$$
f_{4}\big(z,t\big)=z^{2}g_{2}\big(z,t\big)
$$ for some
$g_{2}(z,t)\in\mathbb{C}[z,t]$. Then $v_{3}(0,t)\ne 0$.
\end{corollary}

\begin{lemma}
\label{lemma:very-bad-points-f4} Suppose that $f_{2}(z,t)=zt$.
Then $f_{4}(0,t)=f_{4}(z,0)=0$.
\end{lemma}

\begin{proof}
We may assume that $f_{4}(z,0)\ne 0$. Let $\mathcal{H}$ be
the~linear system on $X$ that is cut out~by
$$
\lambda x+\mu y+\nu t=0,
$$
where $(\lambda:\mu:\nu)\in\mathbb{P}^{2}$. Then the~base locus of
$\mathcal{H}$ consists of the~point $P$.

Let $\mathcal{R}$ be a~proper transform of $\mathcal{H}$ on
the~threefold $V$. Then the~base locus of $\mathcal{R}$ consists
of a~single point that is not contained in any of the~curves
$\bar{L}_{1},\ldots,\bar{L}_{r}$.

The linear system $\mathcal{R}\vert_{B}$ has not base points,
where $B$ is a~general surface in $\mathcal{B}$. But
$$
R\cdot R\cdot B=2n>0,
$$
where $R$ is a~general surface in $\mathcal{R}$. Then
$\mathcal{R}\vert_{B}$ is not composed from a~pencil, which
implies that the~curve $R\cdot B$ is irreducible and reduced by
the~Bertini theorem.

Let $H$ and $M$ be general surfaces in $\mathcal{H}$ and
$\mathcal{M}$, respectively. Then $M\cdot H$ is irreducible and
reduced. Thus, the~linear system $\mathcal{M}|_{H}$ is a~pencil.

The surface $H$ contains no lines passing through $P$, and $H$
can~be given by
$$
w^{3}x+w^{2}y^{2}+w\Big(y^{2}l_{1}\big(x,y,z\big)+xl_{2}\big(x,y,z\big)\Big)+l_{4}\big(x,y,z\big)=0\subset\mathrm{Proj}\Big(\mathbb{C}[x,y,z,w]\Big)\cong\mathbb{P}^{3},
$$
where $l_{i}(x,y,z)$ is a~homogeneous polynomials of degree $i$.

Arguing as in Example~\ref{example:Iskovskikh}, we see that there
is a~pencil $\mathcal{Q}$ on the~surface $H$ such that
$$
\mathcal{Q}\sim\mathcal{O}_{\mathbb{P}^{3}}\big(2\big)\Big\vert_{H},
$$
general curve in $\mathcal{Q}$ is irreducible, and
$\mathrm{mult}_{P}(\mathcal{Q})=4$. Arguing as in the~proof of
Lemma~\ref{lemma:4n-square}, we see that
$\mathcal{M}\vert_{H}=\mathcal{Q}$ by \cite[Theorem~2.2]{ChPa05h}.
Let $M$ be a~general surface in $\mathcal{M}$. Then
$$
M\equiv -2K_{X},
$$
and  $\mathrm{mult}_{P}(M)=4$. The surface $M$ is cut out on $X$
by an~equation
$$
\lambda
x^{2}+x\Big(A_{0}+A_{1}\big(y,z,t\big)\Big)+B_{2}\big(y,z,t\big)+B_{1}\big(y,z,t\big)+B_{0}=0,%
$$
where $A_{i}$ and $B_{i}$ are homogeneous polynomials of degree
$i$, and $\lambda\in\mathbb{C}$.

It follows from $\mathrm{mult}_{P}(M)=4$ that
$B_{1}(y,z,t)=B_{0}=0$.

The coordinated $(y,z,t)$ are also local coordinates on $X$ near
the~point $P$. Then
$$
x=-y^{2}-y\Big(zt+yp_{1}\big(y,z,t\big)\Big)+\text{higher order terms},%
$$
which is a~Taylor power series for $x=x(y,z,t)$, where
$p_{1}(y,z,t)$ is a~linear form.

The surface $M$ is locally given by the~analytic equation
$$
\lambda y^{4}+\Big(-y^{2}-yzt-y^{2}p_{1}\big(y,z,t\big)\Big)\Big(A_{0}+A_{1}\big(y,z,t\big)\Big)+B_{2}\big(y,z,t\big)+\text{higher order terms}=0,%
$$
and $\mathrm{mult}_{P}(M)=4$. Hence, we see that
$B_{2}(y,z,t)=A_{0}y^{2}$ and
$$
A_{1}\big(y,z,t\big)y^{2}+A_{0}y\Big(zt+yp_{1}\big(y,z,t\big)\Big)=0,
$$
which implies that $A_{0}=A_{1}(y,z,t)=B_{2}(y,z,t)=0$. Hence, we
see that a~general~surface in the~pencil $\mathcal{M}$ is cut out
on $X$ by the~equation $x^{2}=0$, which is a~absurd.
\end{proof}

Arguing as in the~proof of Lemma~\ref{lemma:very-bad-points-f4},
we obtain the~following corollary.

\begin{corollary}
\label{corollary:very-bad-points-f4} Suppose that
$f_{2}(z,t)=z^{2}$. Then $f_{4}(0,t)=0$.
\end{corollary}

Let $\mathcal{R}$ be the~linear system on the~threefold $X$ that
is cut out by cubics
$$
xh_{2}\big(x,y,z,t\big)+\lambda \Big(w^{2}x+wy^{2}+q_{3}\big(x,y,z,t\big)\Big)=0,%
$$
where $h_{2}$ is a~form of degree $2$, and $\lambda\in\mathbb{C}$.
Then $\mathcal{R}$ has no fixed components.

Let $M$ and $R$ be general surfaces in $\mathcal{M}$ and
$\mathcal{R}$, respectively. Put
$$
M\cdot R=\sum_{i=1}^{r}m_{i}L_{i}+\Delta,
$$
where $m_{i}\in\mathbb{N}$, and $\Delta$ is a~cycle, whose support
contains no lines among $L_{1},\ldots,L_{r}$.

\begin{lemma}
\label{lemma:Delta-not-empty} The cycle $\Delta$ is not trivial.
\end{lemma}

\begin{proof}
Suppose that $\Delta=0$. Then $\mathcal{M}=\mathcal{R}$ by
\cite[Theorem~2.2]{ChPa05h}.  But $\mathcal{R}$ is not a~pencil.
\end{proof}

We have $\mathrm{mult}_{P}(\Delta)\geqslant
8n-\sum_{i=1}^{r}m_{i}$, because
$\mathrm{mult}_{P}(\mathcal{M})=2n$ and
$\mathrm{mult}_{P}(\mathcal{R})\geqslant 4$. Then
$$
\mathrm{deg}\big(\Delta\big)=12n-\sum_{i=1}^{r}m_{i}\geqslant 2\mathrm{mult}_{P}\big(\Delta\big)\geqslant 2\Big(8n-\sum_{i=1}^{r}m_{i}\Big)%
$$
by Lemma~\ref{lemma:curves-2mult-deg}, because
$\mathrm{Supp}(\Delta)$ does not contain any of the~lines
$L_{1},\ldots,L_{r}$.

\begin{corollary}
\label{corollary:sum-m-i-4n} The inequality
$\sum_{i=1}^{r}m_{i}\geqslant 4n$ holds.
\end{corollary}

Let $R_{1}$ and $R_{2}$ be general surfaces in the~linear system
$\mathcal{R}$. Then
$$
m_{i}\leqslant \mathrm{mult}_{L_{i}}\Big(R_{1}\cdot R_{2}\Big)\mathrm{mult}_{L_{i}}(M)\leqslant\mathrm{mult}_{L_{i}}\Big(R_{1}\cdot R_{2}\Big)n\big\slash 2%
$$
for every $1\leqslant i\leqslant 4$ by
Lemmas~\ref{lemma:cornerstone} and~\ref{lemma:lines-inequalities}.
Then
$$
4n\leqslant\sum_{i=1}^{r}m_{i}\leqslant\sum_{i=1}^{r}\mathrm{mult}_{L_{i}}\Big(R_{1}\cdot R_{2}\Big)n\big\slash 2.%
$$

\begin{corollary}
\label{corollary:sum-m-i-8} The inequality
$\sum_{i=1}^{r}\mathrm{mult}_{L_{i}}(R_{1}\cdot R_{2})\geqslant 8$
holds.
\end{corollary}

Now we suppose that $R_{1}$ is cut out on the~quartic $X$ by the~
equation
$$
w^{2}x+wy^{2}+q_{3}\big(x,y,z,t\big)=0,
$$
and $R_{2}$ is cut out by $xh_{2}\big(x,y,z,t\big)=0$, where
$h_{2}$ is sufficiently general. Then
$$
\sum_{i=1}^{r}\mathrm{mult}_{L_{i}}\Big(R_{1}\cdot T\Big)=\sum_{i=1}^{r}\mathrm{mult}_{L_{i}}\Big(R_{1}\cdot R_{2}\Big)\geqslant 8,%
$$
where $T$ is the~hyperplane section of the~hypersurface $X$ that
is cut out by $x=0$. But
$$
R_{1}\cdot T=Z_{1}+Z_{2},
$$
where $Z_{1}$ and $Z_{2}$ are cycles on $X$ such that $Z_{1}$ is
cut out by $x=y=0$, and $Z_{2}$ is cut~out~by
$$
x=wy+f_{2}\big(z,t\big)+yc_{1}\big(x,y,z,t\big)=0.
$$

\begin{lemma}
\label{lemma:cycle-Z1} The equality
$\sum_{i=1}^{r}\mathrm{mult}_{L_{i}}(Z_{1})=4$ holds.
\end{lemma}

\begin{proof}
The lines $L_{1},\ldots,L_{r}\subset\mathbb{P}^{4}$ are given by
the~equations
$$
x=y=q_{4}\big(x,y,z,t\big)=0,
$$
which implies that $\sum_{i=1}^{r}\mathrm{mult}_{L_{i}}(Z_{1})=4$.
\end{proof}

Hence, we see that
$\sum_{i=1}^{r}\mathrm{mult}_{L_{i}}(Z_{2})\geqslant 4$. But
$Z_{2}$ can~be considered as a~cycle
$$
wy+f_{2}\big(z,t\big)+yc_{1}\big(y,z,t\big)=f_{4}\big(z,t\big)+yv_{3}\big(z,t\big)+y^{2}c_{2}(y,z,t)=0\subset\mathrm{Proj}\Big(\mathbb{C}[y,z,t,w]\Big)\cong\mathbb{P}^{3},%
$$
and, putting $u=w+c_{1}(y,z,t)$, we see that $Z_{2}$ can~be
considered as a~cycle
$$
uy+f_{2}\big(z,t\big)=f_{4}\big(z,t\big)+yv_{3}\big(z,t\big)+y^{2}c_{2}(y,z,t)=0\subset\mathrm{Proj}\Big(\mathbb{C}[y,z,t,u]\Big)\cong\mathbb{P}^{3},%
$$
and we can~consider the~set of lines $L_{1},\ldots,L_{r}$ as
the~set in $\mathbb{P}^{3}$ given by $y=f_{4}(z,t)=0$.

\begin{lemma}
\label{lemma:zt} The inequality $f_{2}(z,t)\ne zt$ holds.
\end{lemma}

\begin{proof}
Suppose that $f_{2}(z,t)=zt$. Then it follows from
Lemma~\ref{lemma:very-bad-points-f4} that
$$
f_{4}\big(z,t\big)=zt\big(\alpha_{1}z+\beta_{1}t\big)\big(\alpha_{2}z+\beta_{2}t\big)
$$
for some
$(\alpha_{1}:\beta_{1})\in\mathbb{P}^{1}\ni(\alpha_{2}:\beta_{2})$.
Then $Z_{2}$ can~be given by
$$
uy+zt=yv_{3}\big(z,t\big)+y^{2}c_{2}\big(y,z,t\big)-uy\big(\alpha_{1}z+\beta_{1}t\big)\big(\alpha_{2}z+\beta_{2}t\big)=0\subset\mathrm{Proj}\Big(\mathbb{C}[y,z,t,u]\Big)\cong\mathbb{P}^{3},%
$$
which implies $Z_{2}=Z_{2}^{1}+Z_{2}^{2}$, where $Z_{2}^{1}$ and
$Z_{2}^{2}$ are cycles in $\mathbb{P}^{3}$ such that $Z_{2}^{1}$
is given~by
$$
y=uy+zt=0,
$$
and $Z_{2}^{2}$ is given by
$uy+zt=v_{3}(z,t)+yc_{2}(y,z,t)-u(\alpha_{1}z+\beta_{1}t)(\alpha_{2}z+\beta_{2}t)=0$.

We may assume that $L_{1}$ is given by $y=z=0$, and $L_{2}$ is
given by $y=t=0$. Then
$$
Z_{2}^{1}=L_{1}+L_{2},
$$
which implies that
$\sum_{i=1}^{r}\mathrm{mult}_{L_{i}}(Z_{2}^{2})\geqslant 2$.

Suppose that $r=4$. Then $\alpha_{1}\ne 0$, $\beta_{1}\ne 0$,
$\alpha_{2}\ne 0$, $\beta_{2}\ne 0$. Hence, we see that
$$
L_{1}\not\subseteq\mathrm{Supp}\big(Z_{2}^{2}\big)\not\supseteq L_{2},%
$$
because
$v_{3}(z,t)+yc_{2}(y,z,t)-u(\alpha_{1}z+\beta_{1}t)(\alpha_{2}z+\beta_{2}t)$
does not vanish on $L_{1}$ and $L_{2}$. But
$$
L_{3}\not\subseteq\mathrm{Supp}\big(Z_{2}^{2}\big)\not\supseteq L_{4},%
$$
because $zt$ does not vanish on $L_{3}$ and $L_{4}$. Then
$\sum_{i=1}^{r}\mathrm{mult}_{L_{i}}(Z_{2}^{2})=0$, which is
impossible.

Suppose that $r=3$. We may assume that
$(\alpha_{1},\beta_{1})=(1,0)$, but $\alpha_{2}\ne 0\ne
\beta_{2}$. Then
$$
L_{2}\not\subseteq\mathrm{Supp}\big(Z_{2}^{2}\big),%
$$
because $v_{3}(z,t)+yc_{2}(y,z,t)-uz(\alpha_{2}z+\beta_{2}t)$ does
not vanish on $L_{2}$. We have
$$
f_{4}\big(z,t\big)=z^{2}t(\alpha_{2}z+\beta_{2}t),
$$
which implies that $v_{3}(0,t)\ne 0$ by
Corollary~\ref{corollary:very-bad-points-derivatives}. Hence, wee
see that
$$
L_{1}\not\subseteq\mathrm{Supp}\big(Z_{2}^{2}\big)\not\supseteq L_{3},%
$$
because $v_{3}(z,t)+yc_{2}(y,z,t)-uz(\alpha_{2}z+\beta_{2}t)$ and
$zt$ do not vanish on $L_{1}$ and $L_{3}$, respectively, which
implies that $\sum_{i=1}^{r}\mathrm{mult}_{L_{i}}(Z_{2}^{2})=0$.
The latter is a~contradiction.

We see that $r=2$. We may assume that
$(\alpha_{1},\beta_{1})=(1,0)$, and either $\alpha_{2}=0$ or
$\beta_{2}=0$.

Suppose that $\alpha_{2}=0$. Then
$f_{4}(z,t)=\beta_{2}z^{2}t^{2}$. By
Lemma~\ref{lemma:very-bad-points-derivatives} and
Corollary~\ref{corollary:very-bad-points-derivatives}, we get
$$
v_{3}\big(0,t\big)\ne 0\ne v_{3}\big(z,0\big),
$$
which implies that $v_{3}(z,t)+yc_{2}(y,z,t)-\beta_{2}zt$ does not
vanish on neither $L_{1}$ nor $L_{2}$. Then
$$
L_{1}\not\subseteq\mathrm{Supp}\big(Z_{2}^{2}\big)\not\supseteq L_{2},%
$$
which implies that
$\sum_{i=1}^{r}\mathrm{mult}_{L_{i}}\big(Z_{2}^{2}\big)=0$, which
is a~contradiction.

We see that $\alpha_{2}\ne 0$ and $\beta_{2}=0$. We have
$f_{4}(z,t)=\alpha_{2}z^{3}t$. Then
$$
v_{3}\big(0,t\big)\ne 0
$$
by Corollary~\ref{corollary:very-bad-points-derivatives}. Then
$L_{1}\not\subseteq\mathrm{Supp}(Z_{2}^{2})$ because
the~polynomial
$$
v_{3}(z,t)+yc_{2}(y,z,t)-\alpha_{2}z^{2}
$$
does not
vanish~on~$L_{1}$.

The line $L_{2}$ is given by the~equations $y=t=0$. But $Z_{2}$ is
given by the~equations
$$
uy+zt=v_{3}\big(z,t\big)+yc_{2}\big(y,z,t\big)-\alpha_{2}uz^{2}=0,%
$$
which implies that $L_{2}\not\subseteq\mathrm{Supp}(Z_{2}^{2})$.
Then $\sum_{i=1}^{r}\mathrm{mult}_{L_{i}}(Z_{2}^{2})=0$, which is
a~contradiction.
\end{proof}

Therefore, we see that $f_{2}(z,t)=z^{2}$. It follows from
Corollary~\ref{corollary:very-bad-points-f4} that
$$
f_{4}\big(z,t\big)=zg_{3}\big(z,t\big)
$$
for some $g_{3}(z,t)\in\mathbb{C}[z,t]$. We may assume that
$L_{1}$ is given by $y=z=0$.

\begin{lemma}
\label{lemma:z2-g3} The equality $g_{3}(0,t)=0$ holds.
\end{lemma}

\begin{proof}
Suppose that  $g_{3}(0,t)\ne 0$. Then
$\mathrm{Supp}(Z_{2})=L_{1}$, because $Z_{2}$ is given by
$$
uy+z^{2}=zg_{3}(z,t)+yv_{3}\big(z,t\big)+y^{2}c_{2}(y,z,t)=0,%
$$
and the~lines $L_{2},\ldots,L_{r}$ are given by the~equations
$y=g_{3}(z,t)=0$.

The cycle $Z_{2}+L_{1}$ is given by the~equations
$$
uy+z^{2}=z^{2}g_{3}\big(z,t\big)+zyv_{3}\big(z,t\big)+zy^{2}c_{2}\big(y,z,t\big)=0,%
$$
which implies that the~cycle $Z_{2}+L_{1}$ can~be given by the~
equations
$$
uy+z^{2}=zyv_{3}\big(z,t\big)+zy^{2}c_{2}\big(y,z,t\big)-uyg_{3}\big(z,t\big)=0.%
$$

We have $Z_{2}+L_{1}=C_{1}+C_{2}$, where $C_{1}$ and $C_{2}$ are
cycles in $\mathbb{P}^{3}$ such that $C_{1}$ is given by
$$
y=uy+z^{2}=0,
$$
and the~cycle $C_{2}$ is given by the~equations
$$
uy+z^{2}=zv_{3}\big(z,t\big)+zyc_{2}\big(y,z,t\big)-ug_{3}\big(z,t\big)=0.%
$$

We have $C_{1}=2L_{2}$. But
$L_{1}\not\subseteq\mathrm{Supp}(C_{2})$ because the~polynomial
$$
zv_{3}\big(z,t\big)+zyc_{2}\big(y,z,t\big)-ug_{3}\big(z,t\big)
$$
does not vanish on $L_{1}$, because $g_{3}(0,t)\ne 0$. Then
$$
Z_{2}+L_{1}=2L_{2},
$$
which implies that $Z_{2}=L_{1}$. Then
$\sum_{i=1}^{r}\mathrm{mult}_{L_{i}}(Z_{2})=1$, which is a~
contradiction.
\end{proof}

Thus, we see that $r\leqslant 3$ and
$$
f_{4}\big(z,t\big)=z^{2}\big(\alpha_{1}z+\beta_{1}t\big)\big(\alpha_{2}z+\beta_{2}t\big)
$$
for some
$(\alpha_{1}:\beta_{1})\in\mathbb{P}^{1}\ni(\alpha_{2}:\beta_{2})$.
Then
$$
v_{3}\big(0,t\big)\ne 0
$$
by Corollary~\ref{corollary:very-bad-points-derivatives}. But
$Z_{2}$ can~be given~by~the~equations
$$
uy+z^{2}=yv_{3}\big(z,t\big)+y^{2}c_{2}\big(y,z,t\big)-uy\big(\alpha_{1}z+\beta_{1}t\big)\big(\alpha_{2}z+\beta_{2}t\big)=0\subset\mathrm{Proj}\Big(\mathbb{C}[y,z,t,u]\Big)\cong\mathbb{P}^{3},%
$$
which implies $Z_{2}=Z_{2}^{1}+Z_{2}^{2}$, where $Z_{2}^{1}$ and
$Z_{2}^{2}$ are cycles on $\mathbb{P}^{3}$ such that $Z_{2}^{1}$
is given by
$$
y=uy+z^{2}=0,
$$
and the~cycle $Z_{2}^{2}$ is given by the~equations
$$
uy+z^{2}=v_{3}\big(z,t\big)+yc_{2}\big(y,z,t\big)-u\big(\alpha_{1}z+\beta_{1}t\big)\big(\alpha_{2}z+\beta_{2}t\big)=0,%
$$
which implies that $Z_{2}^{1}=2L_{1}$. Thus, we see that
$\sum_{i=1}^{r}\mathrm{mult}_{L_{i}}(Z_{2}^{2})\geqslant 2$.

\begin{lemma}
\label{lemma:z2-r-3} The inequality $r\ne 3$ holds.
\end{lemma}

\begin{proof}
Suppose that $r=3$. Then $\beta_{1}\ne 0\ne\beta_{2}$, which
implies that
$$
L_{1}\not\subseteq\mathrm{Supp}\big(Z_{2}^{2}\big),%
$$
because
$v_{3}(z,t)+yc_{2}(y,z,t)-u(\alpha_{1}z+\beta_{1}t)(\alpha_{2}z+\beta_{2}t)$
does not vanish on $L_{1}$. But
$$
L_{2}\not\subseteq\mathrm{Supp}\big(Z_{2}^{2}\big)\not\supseteq L_{3},%
$$
because $\beta_{1}\ne 0\ne\beta_{2}$. Then
$\sum_{i=1}^{r}\mathrm{mult}_{L_{i}}(Z_{2}^{2})=0$, which is a~
contradiction.
\end{proof}

Thus, we see that either $r=1$ or $r=2$.

\begin{lemma}
\label{lemma:z2-r-2} The inequality $r\ne 2$ holds.
\end{lemma}

\begin{proof}
Suppose that $r=2$. We may assume that
\begin{itemize}
\item  either $\beta_{1}\ne 0=\beta_{2}$,%
\item or $\alpha_{1}=\alpha_{2}$ and $\beta_{1}=\beta_{2}\ne 0$.%
\end{itemize}

Suppose that $\beta_{2}=0$. Then
$f_{4}(z,t)=\alpha_{2}z^{3}(\alpha_{1}z+\beta_{1}t)$ and
$$
L_{1}\not\subseteq\mathrm{Supp}\big(Z_{2}^{2}\big),%
$$
because
$v_{3}(z,t)+yc_{2}(y,z,t)-\alpha_{2}uz(\alpha_{1}z+\beta_{2}t)$
does not vanish on $L_{1}$. But $L_{2}$ is given by
$$
y=\alpha_{1}z+\beta_{1}t=0,
$$
which implies that $z^{2}$ does not vanish on $L_{2}$, because
$\beta_{1}\ne 0$. Then
$$
L_{2}\not\subseteq\mathrm{Supp}\big(Z_{2}^{2}\big),
$$
which~implies that
$\sum_{i=1}^{r}\mathrm{mult}_{L_{i}}(Z_{2}^{2})=0$, which is
a~contradiction.

Hence, we see that $\alpha_{1}=\alpha_{2}$ and
$\beta_{1}=\beta_{2}\ne 0$. Then
$L_{1}\not\subseteq\mathrm{Supp}(Z_{2}^{2})$, because
$$
v_{3}(z,t)+yc_{2}(y,z,t)-u(\alpha_{1}z+\beta_{1}t)^{2}
$$ does not
vanish~on~$L_{1}$. But
$L_{2}\not\subseteq\mathrm{Supp}(Z_{2}^{2})$, because $z^{2}$ does
not vanish on $L_{2}$. Then
$$
\sum_{i=1}^{r}\mathrm{mult}_{L_{i}}\big(Z_{2}^{2}\big)=0,
$$
which is a~contradiction.
\end{proof}

We see that $f_{4}(z,t)=z^{2}$ and $f_{4}(z,t)=\mu z^{4}$ for some
$0\ne \mu\in\mathbb{C}$. Then $Z_{2}^{2}$ is given by
$$
uy+z^{2}=v_{3}\big(z,t\big)+yc_{2}\big(y,z,t\big)-\mu z^{2}=0,%
$$
where $v_{3}(0,t)\ne 0$ by
Corollary~\ref{corollary:very-bad-points-derivatives}. Thus, we
see that $L_{1}\not\subseteq\mathrm{Supp}(Z_{2}^{2})$, because
$$
v_{3}\big(z,t\big)+yc_{2}\big(y,z,t\big)-\mu z^{2}
$$ does not
vanish on $L_{1}$. Then
$\sum_{i=1}^{r}\mathrm{mult}_{L_{i}}(Z_{2}^{2})=0$, which is
a~contradiction.

The assertion of Proposition~\ref{proposition:very-bad-points} is
proved.

The assertion of Theorem~\ref{theorem:main} follows from
Propositions~\ref{proposition:curves},
\ref{proposition:good-points}, \ref{proposition:bad-points},
\ref{proposition:very-bad-points}.

\end{document}